\documentclass[a4paper,11pt,fullpage]{article}
\usepackage[english]{babel}

\usepackage[left=1in,right=1in,top=1in,bottom=1in]{geometry}

\usepackage{colortbl}
\usepackage{amsmath, amsthm, amscd, amsfonts, amssymb, graphicx, color}
\usepackage{indentfirst}
\usepackage{amsmath}
\usepackage{amsthm}
\usepackage{amsfonts}
\usepackage{amssymb}
\usepackage{amsxtra}
\usepackage{latexsym}
\usepackage{ifthen}
\usepackage{cite}
\usepackage{esint}
\usepackage[cp1250]{inputenc}

\usepackage{epic}
\usepackage{eepic}

\numberwithin{equation}{section}
\newtheorem{theorem}{Theorem}[section]

\newtheorem{remark}{Remark}[section]
\newtheorem{definition}{Definition}[section]

\def\XXint#1#2#3{{\setbox0=\hbox{$#1{#2#3}{\int}$}
     \vcenter{\hbox{$#2#3$}}\kern-.5\wd0}}

\begin{document}

\title{Potential analysis of multi layers optic fiber models}

\author{Kateryna Buryachenko\\
 Humboldt-Universit\"at zu Berlin, Berlin, Germany,\\
Vasyl' Stus Donetsk National University, Vinnytsia, Ukraine.\\
Yuliya Kudrych\\
Vasyl' Stus Donetsk National University, Vinnytsia, Ukraine}

%\shorttitle{Short paper title for the headers}

%\shortauthor{F. Author, S. Author}

\date{}

\maketitle

\begin{abstract}
This work is devoted to the development of qualitative methods for the study of nonlinear heterogeneous structures, models of which are elliptic equations, which describe complex nonlinear processes in heterogeneous media. They may also include the structures, consisting of several parts (phases or layers): multiphase solid and liquid materials; optic fiber and optic cable layers, anisotropic medium, etc. Relevance of the chosen direction is due to the fact that many processes in heterogeneous environments under conditions of high temperatures, heavy loads and significant deformations are described using nonlinear differential equations with discontinuous (singular) data (coefficients, right-hand side, boundary and initial conditions, etc.). At the same time, the concept of weak solutions that meet the modern needs of mathematical physics arose. Nonlinear differential equations have a complex structure, which actually makes them impossible to study by finding solutions in an explicit form. Therefore, the development of qualitative  methods for their investigations becomes an extremely important tool. This paper considers mathematical models of multilayer optic fiber and cable, which consist of 3 and 5 different materials respectively with different properties. Using potential theory, the behavior of a weak solution of this equation at a fixed point is estimated and analyzed by the value of the nonlinear Wolff potential from the right hand side. We study pointwise properties that play a key role in the further study: expansion of positivity Harnack's inequalities, regularities and others. The paper discusses also the application of the obtained theoretical results for the problem of modeling and analyzing of  optic fiber and optic cable modern technologies.

\end{abstract}

Laplace, Wolff potential, weak solution, pointwise estimates.

\bigskip
\section{Introduction}

We focus here on the development of qualitative methods of nonlinear analysis for the study of  double-phase elliptic equations with variable exponents and their applications in modern optic technologies. The active development of the problem under consideration is evidenced by numerous high citing publications during the past 2-3 years in leading journals: V. B\"ogelein, F. Duzaar, P. Marcellini, C. Scheven \cite{BDM}, V. B\"ogelein,  M.Strunk \cite{BS}, C. De Filippis, G. Mingione \cite{FM} and others.  The double-phase elliptic equations of the divergence
form were studied in first in the papers \cite{Sb_37, Sb_38} as
models of strictly anisotropic materials and for the description
of Lavrent'ev phenomenon. H\"{o}lder continuity and Harnack's
inequality for bounded solutions to the homogeneous equation  were obtained in \cite{Sb_17},
\cite{Sb_19} under the same conditions, which we have herein .

These works were fundamental for further studies of the existence and regularity of solutions of various types of problems for such equations. The novelty of the results of this paper is the development of new functional methods of nonlinear analysis for the study of new actual problems, the mathematical models of which are double-phase elliptic equations with variable exponents. We introduce here the new potential estimates for the weak nonnegative solutions via nonlinear Wolff potential of the right hand side $f\in L^1$ of the equation and discuss their applications in the modeling of optic fiber devices.

The considering class of double-phase equations serves as mathematical model of media including structures which consist of several parts (phases or layers): multiphase solid and liquid materials; porous, anisotropic media; optic fiber layers, optic cable layers, light diodes, semiconductors devices, etc. The relevance of the chosen direction is due to the fact that many processes in heterogeneous environments under conditions of high temperature, heavy loads and significant deformations are described by using similar equations and with discontinuous (singular) data (coefficients, right-hand side, boundary and initial conditions, etc.). In our case this is a right hand side $f\in L^1.$ At the same time, the concept of weak solutions is widely used, which meets the modern needs of mathematical physics. Nonlinear differential equations have a complex structure, which actually makes it impossible to study them by finding solutions in an explicit form. Therefore, the development of qualitative methods of analysis becomes an extremely important tool. In the present manuscript we obtain new pointwise estimates for the weak nonnegative solution via nonlinear Wolff potential from the right hand side of elliptic equations with non-standard growth conditions, $(p, q)-$ double-phase equations, with variable exponents: $p(x), q (x)$. Obtained in the manuscript new pointwise properties for the weak solutions via nonlinear Wolff potential from the right-hand side $f\in L^1$ will explore fundamental qualitative properties that play a key role in further studying the behavior of solutions: boundedness, expansion of positivity, H\"older continuity, and Harnack's inequalities. 

The main results of the current paper are  expansions of the works \cite{BS} and \cite{KS} for the case of double-phase elliptic equations with variable exponents $p(x),\, q(x).$

We consider also  mathematical models of multilayer optic fiber and multilayer optic cable, which consist of 3 and 5 different materials with different properties and discuss  the application of the obtained theoretical results for the problem of modeling and analyzing of optic fiber and optic cable modern technologies.

\section{Statement of the problem}
In a bounded domain $\Omega\subset \mathbb{R}^n,\,n\geq 2$ we
consider a double-phase  elliptic equation with variable exponents:
\begin{equation}\label{eq2.0}
-{\rm div} \left[(|\nabla u|^{p(x)-2}+a(x)|\nabla u|^{q(x)-2})\nabla u\right]=f(x)\geq 0,
\end{equation}

\begin{equation}\label{eq2.2}
-{\rm div} A(x,\,\nabla u)=f(x)\geq 0,
\end{equation}
where $f(x)\in L^1(\Omega)$. We assume that  the function
$A(x,\,\xi)=|\xi|^{p(x)-1}+a(x)|\xi|^{q(x)-1}:\Omega \times \mathbb{R}^n\to  \mathbb{R}^n $
satisfies the conditions

1) $A(x,\,\xi)$ satisfies the Carath\'{e}odory condition,

2) $A(x,\xi)\xi\geq \mu_1 (|\xi|^{p(x)}+a(x)|\xi|^{q(x)}),$

3) $|A(x,\xi)|\leq \mu_2 (|\xi|^{p(x)-1}+a(x)|\xi|^{q(x)-1}),$

with some constants $\mu_1,\,\mu_2>0.$ 

We also assume that
$$0\leq a(x)\in C^{0,\,\alpha}(\Omega),\quad \alpha\in
(0,\,1].$$

Let $\cal M$ be a set of all measurable functions, $p(x), q(x):\Omega\to (1,\infty).$ For $p(x), q(x)\in \cal M,$ we set:

$$p_-:={\rm ess inf}_{x\in\Omega} p(x),\,q_-= {\rm ess inf}_{x\in\Omega} q(x), p_+:={\rm ess sup}_{x\in\Omega} p(x),\,q_+= {\rm ess sup}_{x\in\Omega} q(x).$$

We assume the following for the powers of nonlinearity:
\begin{equation}\label{cond_powers}
1<p_-\leq p_+\leq q_-\leq q_+\leq \min\left(p_-+\alpha,\,
\frac{n(p_--1)}{n-p_-}\right),\quad q_+<n.
\end{equation}

Let us introduce the necessary definitions.

\begin{definition}\label{de1} Let $G(x,t)=t(t^{p(x)-1}+a(x)t^{q(x)-1})$. Then $W^{1,G}(\Omega)$ denotes the class of functions $u$ that are weakly differentiable
in $\Omega$ and satisfy the condition
$$ \int\limits_{\Omega}G(a(x),\,|\nabla u|)\,dx<\infty.$$
\end{definition}

\vskip3.5mm
\begin{definition}\label{de2}
We say that $u$ is a weak solution to Eq. ($\ref{eq2.2}$), if $u\in
W^{1,G}(\Omega)$  and it satisfies the integral identity
\begin{equation}\label{eq2.5}
\int\limits_{\Omega}A(x,\,\nabla
u)\nabla\varphi\,dx=\int\limits_{\Omega}f\,\varphi\,dx,
\end{equation}
 for all
$\varphi \in \overset {0}W^{1,G}(\Omega)$.
\end{definition}

In the case of Eq.(\ref{eq2.0}) condition (\ref{eq2.5}) takes the form:

\begin{equation}\label{de3}
\int\limits_{\Omega}\left(|\nabla u|^{p(x)-1}+a(x)|\nabla u|^{q(x)-1}\right)\nabla\varphi\,dx=\int\limits_{\Omega}f\,\varphi\,dx.
\end{equation}
\vskip10mm
We will prove the pointwise estimates for a nonnegative weak solution
to the double-phase equation (\ref{eq2.0}) in terms of the nonlinear Wolff
potentials:

\begin{equation*}
W_{1,p(x)}^f(x_0,\,R)=\sum_{j=0}^{\infty}\left(\rho_j^{p(x)-n}\int\limits_{B_{\rho_j}(x_0)}f\,dx\right)^{\frac{1}{p(x)-1}},\,\rho_j=\frac{R}{2^j},\,j=0,\,1,\,...
\end{equation*}
\begin{equation*}
W_{1,q(x)}^f(x_0,\,R)=\sum_{j=0}^{\infty}\left(\rho_j^{q(x)-n}\int\limits_{B_{\rho_j}(x_0)}f\,dx\right)^{\frac{1}{q(x)-1}},\,\rho_j=\frac{R}{2^j},\,j=0,\,1,\,...,
\end{equation*}
under assumption that the series in the above formulae are
convergent, i.e. the Wolff potentials are finite.

Let us note that double-phase elliptic equations of the divergence
form were studied in first in the papers \cite{Sb_37, Sb_38} as
models of strictly anisotropic materials and for the description
of Lavrent'ev phenomenon. H\"{o}lder continuity and Harnack
inequality for bounded solutions to the homogeneous equation (\ref{eq2.0})
(with function  $f \equiv 0$)  were obtained in \cite{Sb_17},
\cite{Sb_19} under conditions (\ref{cond_powers}).

\section{Main result}
The main result of the present work is the following theorem.

\begin{theorem}\label{t1}
Let  $u\in W^{1,G}(\Omega)\cap L^{\infty}$ be a nonnegative weak
solution to Eq. (\ref{eq2.0}). Let conditions (\ref{cond_powers}) be satisfied
and let
$[a]_{C^{0,\alpha}(\Omega)}:=\underset{x,y\in\Omega,\,x\neq y}{\rm
sup}\frac{|a(x)-a(y)|}{|x-y|^{\alpha}}.$ Assume also that the
point $x_0\in\Omega$ is such that $B_{4\rho}(x_0)\subset\Omega$.
Then there exist constants $c_1,\,c_2>0$ depending only on
$p_-,\,q_+,\,n,$ $[a]_{C^{0,\alpha}(\Omega)}$ and
$||u||^{q_+-p_-}_{L^{\infty}(\Omega)}$ such that, under condition
$a(x_0)=0$ the following estimate holds:
\begin{equation}\label{eq2.10}
c_1W_{1,p_-}^f(x_0,\rho)\leq u(x_0)\leq
c_2\underset{B_{\rho}(x_0)}{\rm inf}u+c_2
W_{1,p_-}^f(x_0,2\rho).\end{equation}

If $a(x_0)>0$ and  $\rho_0^{\alpha}=\frac{a(x_0)}{4
[a]_{C^{0,\alpha}(\Omega)}}\geq\rho^{\alpha},$ then there exist
constants $c_3,\,c_4>0$ depending on
$p_-,\,q_+,\,n,\,[a]_{C^{0,\alpha}(\Omega)},\,||u||^{q_+-p_-}_{L^{\infty}(\Omega)}$
and $a(x_0)$ such that the following estimate
\begin{equation}\label{eq2.11}
 c_3W_{1,q_+}^f(x_0,\rho)\leq\rho+u(x_0)\leq 3\rho+
c_4\underset{B_{\rho}(x_0)}{\rm
inf}u+c_4W_{1,q_+}^f(x_0,2\rho)\end{equation} holds.

Under conditions $a(x_0)>0$ and $\rho_0<\rho$ will be true the
estimate
\begin{equation*}
c_3W_{1,q_+}^f(x_0,\rho)+c_3(W_{1,p_-}^f(x_0,\rho)-W_{1,p_-}^f(x_0,\rho_0))\leq\rho+u(x_0)\leq
\end{equation*}
\begin{equation}\label{eq2.12}
\leq 3\rho+c_4\underset{B_{\rho}(x_0)}{\rm
inf}u+c_4W_{1,q_+}^f(x_0,2\rho)+c_4(W_{1,p_-}^f(x_0,2\rho)-W_{1,p_-}^f(x_0,2\rho_0)).
\end{equation}
\end{theorem}
 \vskip3.5mm

{\it Proof.}
The result of this theorem will follow from the analogue result, proved in \cite{BS} for the double-phase equation with constant powers of nonlinearity $p, q:$
\begin{equation}\label{eq2.1}
-{\rm div} \left[(|\nabla u|^{p-2}+a(x)|\nabla u|^{q-2})\nabla u\right]=f(x)\geq 0,
\end{equation}
with 

\begin{equation}\label{2}
1<p\leq q\leq \min\left(p+\alpha,\,
\frac{n(p-1)}{n-p}\right),\quad q<n.
\end{equation}

Taking into account our conditions (\ref{cond_powers}), we have the statement of our theorem as a consequence of the analogous result for Eq.(\ref{eq2.1}), see \cite{BS}.

\vspace{2mm}
\begin{remark}\label{r1} In the case  $a(x_0) = 0$ inequality (\ref{eq2.10}) yields the
known result of Kilpel\"ainen and Mal\'y \cite{Sb_1}, where there
were obtained the pointwise estimates of solutions to a
quasilinear elliptic equation with the p-Laplace and  measure
$\mu$ on the right-hand side with the help of the nonlinear Wolff
potential $W_{\beta,\,p_-}^{\mu}(x_0,\,R)$:
\begin{equation}\label{eq2.4}
W_{\beta,\,p_-}^{\mu}(x_0,\,R):=\sum_{j=0}^{\infty}\left(\frac{\mu(B_{\rho_j}(x_0))}{\rho_j^{n-\beta
p_-}}\right)^{\frac{1}{p_--1}},\,\rho_j=\frac{R}{2^j},\,j=0,\,1,\,2,...\end{equation}

\end{remark}

\section{Applications to the multi layers optic fiber models}

Consider the multilayer optic fiber model, described by the exponents:

\begin{equation}\label{p-q}
p(x)=\left\{\begin{array}{cc}
   p_1  & x\in\Omega_1, \\
   p_2  & x\in\Omega_2,\\
   p_3 & x\in\Omega_3,\\
   \cdots &\\
   p_n& x\in \Omega_n;\\
\end{array}\right. \,\,q(x)=\left\{\begin{array}{cc}
   q_1  & x\in\Omega_1, \\
   q_2  & x\in\Omega_2,\\
   q_3 & x\in\Omega_3,\\
   \cdots &\\
   q_n& x\in \Omega_n,\\
\end{array}\right.
\end{equation}
with the constant $p_i,\,q_i,\, i=1,...,n,$ depending on the number of layers in the optic fiber.
In this case of discrete-valued $p(x)$ and $q(x),$ we stand
$$ p_-=\min_{i=1,..,n}{p_i},\,\,\, q_-=\min_{i=1,..,n}{q_i},\,\,\,p_+=\max_{i=1,..,n}{p_i},\,\,\,q_+=\max_{i=1,...,n}{q_i}.$$
As usual,  optic cable is a carefully designed multilayer designed to protect sensitive optic fiber and ensure its optimal performance under various environmental conditions and mechanical loads. The main components working in interaction include:
\begin{itemize}
        \item {\bf Core}: $\Omega_1,$ The innermost part of the cable, serving as a way to transmit light. It is usually made of high-purity glass or, less commonly for single-mode fibers, of plastics;
\item  {\bf Cladding}: $\Omega_2,$ The optical layer immediately surrounding the core. Its material composition is chosen to have a lower refractive index than the core, which is a critical property that contributes to the complete internal reflection and retention of light in the core;
 \item {\bf Buffer layer}: $\Omega_3,$ A protective coating applied directly over the shell. This layer provides substantial physical protection to the fiber, protecting it from minor abrasive damage, impact and exposure to environmental elements;
\item {\bf Power elements}: $\Omega_4,$ These components are strategically integrated into the cable structure to provide tensile strength and mechanical reinforcement, protecting the optic fiber from stretching, bending, and crushing;
\item {\bf Coating}: $\Omega_5,$ The outer protective layer of the cable. This layer provides comprehensive protection against moisture, ultraviolet radiation, chemicals and mechanical damage, and often serves to identify the cable.

\end{itemize}
For example, in the case of a optic fiber it is a carefully designed by three multilayer designed to protect sensitive optic fiber and ensure its optimal performance under various environmental conditions. The main components working in interaction include only two parts (core and cladding). Please, see the following single-mode optic fiber:

\includegraphics[height=6.5cm]{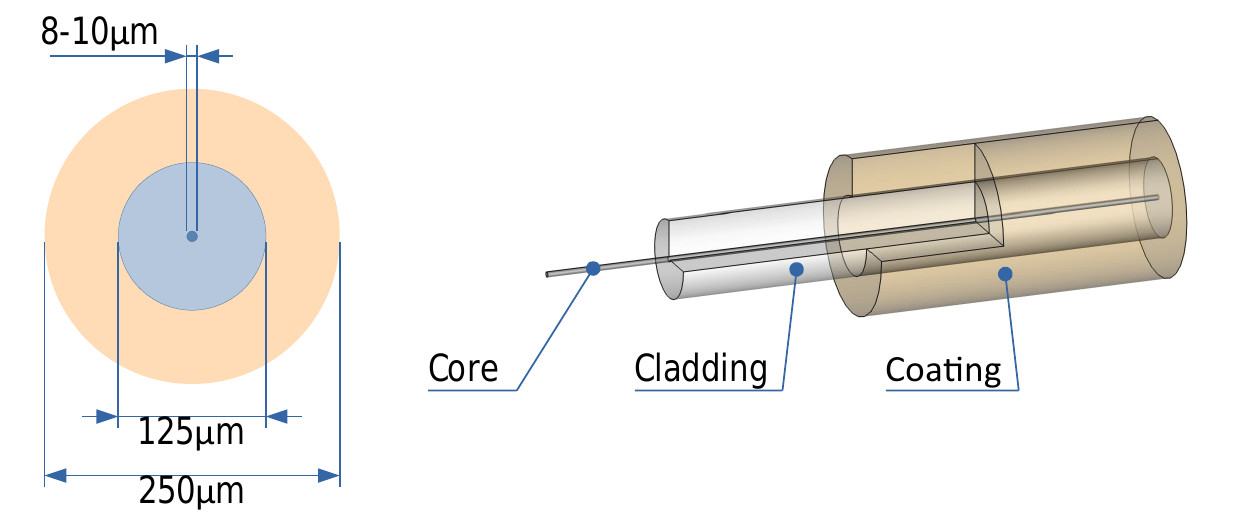}

Taking into account that optic fiber consists of different kinds of materials, it is natural, that the value of powers of nonlinearity, $p(x), q(x)$ take a different values $(p,q)$  on each of layers $\Omega_i,\, i=1,...,5.$ 

Thus, the core usually consists of ultrapure quartz (SiO2). To achieve the required higher refractive index relative to the shell, quartz is precisely doped with elements such as germanium dioxide (GeO2) \cite{FO}. The ultra-purity of quartz glass is paramount to minimize light absorption and scattering, thereby ensuring high transmission efficiency. The cladding layer is usually made of pure quartz or fluorine-doped quartz, which effectively reduces its refractive index compared to the germanium-doped core \cite{FO}.
Acrylate polymers or polyimides are used for buffer. These materials are chosen because of their adhesion to glass and protective properties. Aramid threads (e.g. Kevlar, Twaron) are wide use materials for power elements.The cable outer sheath is the most visible protective layer of the optical cable. Its main role is to protect the internal components from environmental factors, mechanical damage and fire dangers. The choice of sheath material is very application-dependent, balancing performance, cost and safety requirements, for instance: polyvinyl chloride, polyethylene, polyurethane and others \cite{Fiber}.

For (\ref{p-q}) the result of Theorem \ref{t1} can be applied.
So, we can estimate the pointwise value of the solution $u(x_0)$ via nonlinear Wolff potential of the right-hand side $f\in L^1(\Omega)$,  depending on the point $x_0\in\Omega=\Omega_1\cup\Omega_2\cup\Omega_3\cup\Omega_4\cup\Omega_5.$

\vskip5mm

\section*{ Conclusions} The paper discusses a mathematical model of multilayer optic fiber and optic cable, which consists of 3 and 5 different materials respectively with different properties.
Using potential theory, the behavior of a weak solution of this equation at a fixed point from the value of the nonlinear Wolff potential from the right side is analyzed. This result complements the work of one of the author \cite{BS} in the case of variable powers $p(x), q(x)$ of nonlinearity. Additionally, the paper discusses the application of the obtained results for the problem of modeling and analyzing of optic fiber modern technologies.

\vspace{2mm}

{\bf Conflict of interest and ethics.} The authors declare no conflict of interests. The authors also declare  full adherence to all journal research ethics policies, namely involving the participation of human subjects anonymity and consent to publish.

\vspace{2mm}

{\bf Acknowledgements.} 

The authors also thank to their colleagues-physicists Mykola Pasichnyy and Vasyl Komarov for the fruitful discussions and clarifications in area of optic fiber modern technologies.

\section*{Funding}
Yuliya Kudrych is supported within the framework of the program 2025.06 "Science for Strengthening the Defense Capability and National Security of Ukraine" of the National Research Foundation of Ukraine (NRFU), project no. 2025.06/0090, state registration no. 0125U003181.

\end{document}